\documentclass[10pt]{article}

\usepackage[comma]{natbib}
\usepackage{anysize}
\usepackage{xcolor}
\usepackage{amsmath,amssymb,amsthm,amsfonts,amsopn,bm,latexsym,scalerel}
\usepackage{verbatim,graphicx}
\usepackage[latin1]{inputenc} 
\usepackage[english]{babel}

\usepackage{pgfplots}
\pgfplotsset{compat=1.10}
\usepackage{tikz}
\usepackage{tikz-3dplot}
\usetikzlibrary{calc}
\tikzset{viewport/.style 2 args={
		x={({cos(-#1)*1cm},{sin(-#1)*sin(#2)*1cm})},
		y={({-sin(-#1)*1cm},{cos(-#1)*sin(#2)*1cm})},
		z={(0,{cos(#2)*1cm})}
}}

\pgfplotsset{only foreground/.style={
		restrict expr to domain={rawx*\CameraX + rawy*\CameraY + rawz*\CameraZ}{-0.05:100},
}}
\pgfplotsset{only background/.style={
		restrict expr to domain={rawx*\CameraX + rawy*\CameraY + rawz*\CameraZ}{-100:0.05}
}}

\def\addFGBGplot[#1]#2;{
	\addplot3[#1,only background, opacity=0.25] #2;
	\addplot3[#1,only foreground] #2;
}

\newcommand{\ViewAzimuth}{15}
\newcommand{\ViewElevation}{15}
\tdplotsetmaincoords{75}{105} 

\marginsize{3cm}{3cm}{3cm}{3cm}

\newcommand{\bx}{{\bm x}}
\newcommand{\by}{{\bm y}}

\newcommand{\bmu}{{\bm \mu}}

\renewcommand{\ell}{\mathcalboondox{l}}

\newcommand{\bbf}{\mathbf{f}}

\newcommand{\Log}{{\rm Log}}

\theoremstyle{plain}

\title{A note on the geodesic normal distribution\\ on the sphere}

\author{Jos\'e E. Chac\'on\footnote{Departamento de
Matem\'aticas, Universidad de Extremadura, E-06006 Badajoz, Spain. E-mail:
{\tt jechacon@unex.es}} \ and Andrea Meil\'an-Vila\footnote{Departamento de Estad\'\i stica, Universidad Carlos III de Madrid, E-28911 Legan\'es, Spain. Email: {\tt ameilan@est-econ.uc3m.es}}}

\begin{document}

\maketitle

\begin{abstract}
\noindent This paper presents an alternative formulation of the geodesic normal distribution on the sphere, building on the work of \cite{Hauberg2018}. While the isotropic version of this distribution is naturally defined on the sphere, the anisotropic version requires projecting points from the hypersphere onto the tangent space. In contrast, our approach removes the dependence on the tangent space and defines the geodesic normal distribution directly on the sphere. 
Moreover, we demonstrate that the density contours of this distribution are exactly ellipses on the sphere, providing intriguing alternative characterizations for describing this locus of points. 
\end{abstract}


\newpage

\section{Introduction}
\label{sec1}

Directional data refers to observations whose support lies on the unit hypersphere $\mathbb{S}^p:= \{\bx \in \mathbb{R}^{p+1}\colon \|\bx\|^2=\bx^\top \bx = 1\}$ in $\mathbb{R}^{p+1}$, $p \ge 1$.  Isotropic normal distributions are commonly used to model such data, with the resulting density contours forming (hyper)spherical circles on $\mathbb{S}^p$. In this setting, two main geometric frameworks are used to define distances. The extrinsic framework employs the Euclidean norm, where distances between objects are computed after embedding them into Euclidean space \citep{bhattacharya2012nonparametric}. The classical von Mises-Fisher distribution arises from using the Euclidean distance on $\mathbb{S}^p$. In contrast, the intrinsic framework defines distances via the shortest-path geodesics (great circles) on the hypersphere, leading to the (isotropic) geodesic normal distribution on $\mathbb{S}^p$.  This term was first introduced by \cite{coeurjolly2012geodesic} for $p=1$, and later generalized to higher dimensions by \citet{Hauberg2018}, who referred to it as the \emph{spherical normal}. A comprehensive study of this distribution was provided by \cite{you2022parameter}. 

Extensions to anisotropic distributions have been proposed within both frameworks. For example, \citet{Kent1982} introduced a 5-parameter distribution on the sphere $\mathbb{S}^2$ that uses the Euclidean norm, but does not account for the curvature of the sphere. \cite{paine2018elliptically} proposed the elliptically symmetric angular Gaussian model (ESAG), which is computationally more efficient than the Kent model for both simulation and density calculation. Other Kent-like alternatives that offer the same advantages as the ESAG model include the scaled von Mises-Fisher family suggested by \cite{scealy2019scaled}, which incorporates an additional parameter to control tail-weight. Within the intrinsic framework, \cite{Hauberg2018} introduced the anisotropic spherical normal distribution on $\mathbb{S}^p$, which takes into account the curvature of the sphere. This distribution employs the logarithm map to project points from the hypersphere onto its tangent space, building upon the general Riemannian normal distribution described in \citet{pennec2006intrinsic}.

The main aim of this paper is to present an alternative formulation of Hauberg spherical normal distribution that avoids its reliance on the tangent space. In this new formulation, the distribution is expressed directly on the spherical manifold, providing a more natural representation. Furthermore, it is shown that the density contours of this distribution are true ellipses on the sphere, in contrast to other anisotropic distributions, thereby offering a clearer geometric understanding of the distribution.

The rest of this paper is organized as follows. Section \ref{sec:normal} provides a brief overview of isotropic and anisotropic normal distributions on the sphere. In Section \ref{sec:gn}, we present our alternative formulation of the anisotropic geodesic normal distribution. Finally, Section \ref{sec:sph_ell} exhibits several equivalent equations for defining a spherical ellipse and demonstrates that the density contours of the anisotropic geodesic normal distribution correspond to spherical ellipses.

\section{Normal distributions on the sphere}
\label{sec:normal}

For a given distance $d$ on the $p$-dimensional sphere $\mathbb S^p$, the class of isotropic normal distributions with respect to $d$ consists of those whose density is proportional to $\exp\{-\frac{\kappa}2d(\bx,\bmu)^2\}$ with respect to the Lebesgue measure $\sigma_p({\mathsf d}\bx)$ on $\mathbb S^p$. Here, $\bmu\in\mathbb S^p$ represents the mean direction and $\kappa>0$ is a concentration parameter. Their density contours are given by $\{\bx\in\mathbb S^p\colon d(\bx,\bmu)=\rho\}$ for some $\rho>0$. This is precisely the definition, as a locus of points, of a spherical circle with center $\bmu$ and radius $\rho$ (with respect to $d$).

When using the Euclidean distance on the hypersphere, given by $d_E(\bx,\bmu)=\{2(1-\bx^\top\bmu)\}^{1/2}$ for $\bx,\bmu\in\mathbb S^p$, we obtain the von Mises-Fisher distribution, whose density is $f_{\rm vMF}(\bx;\bmu,\kappa)=C_{\rm vMF}(\kappa)^{-1}\exp(\kappa\,\bx^\top\bmu)$. Here, $C_{\rm vMF}(\kappa)=\int_{\mathbb S^p}\exp(\kappa\,\bx^\top\bmu)\sigma_p(\mathsf d\bx)$ denotes the normalizing constant, whose explicit formula is also well known \citep[see][Section 9.3.2]{Mardia1999a}. In contrast, the geodesic distance is defined as $d_G(\bx,\bmu)=\arccos(\bx^\top\bmu)$, which represents the length of the shortest arc between $\bx$ and $\bmu$ along $\mathbb S^p$. Hence, the (isotropic) geodesic normal distribution has density fuction $f_{\rm GN}(\bx;\bmu,\kappa)=C_{\rm GN}(\kappa)^{-1}\exp\{-\frac{\kappa}2\arccos^2(\bx^\top\bmu)\}$. For the case 
$p=1$, an explicit formula for the normalizing constant $C_{\rm GN}(\kappa)=\int_{\mathbb S^p}\exp\{-\frac\kappa2\arccos^2(\bx^\top\bmu)\}\sigma_p(\mathsf d\bx)$ was derived by \cite{coeurjolly2012geodesic}. For higher dimensions, formulas were provided in \citet{Hauberg2018}.

Various approaches have been proposed to extend those existing models to obtain anisotropic normal distributions. \citet{Kent1982} introduced a 5-parameter distribution on the $\mathbb S^2$ with density $f_{\rm K}(\bx;\kappa,\beta,\mathbf \Gamma)=C_{\rm K}(\kappa,\beta)^{-1} \exp\{\kappa\bm\gamma_1^\top\bx+\beta[(\bm\gamma_2^\top\bx)^2-(\bm\gamma_3^\top\bx)^2]\}$, where $\kappa\geq0$, $\beta\geq 0$, $\bm\Gamma=(\bm\gamma_1|\bm\gamma_2|\bm\gamma_3)$ is a $3\times 3$ orthogonal matrix and $C_{\rm K}(\kappa,\beta)$ denotes the normalizing constant, which does not depend on $\bm\Gamma$. This distribution is a generalization of the von Mises-Fisher distribution, which corresponds to the particular choices of $\beta=0$ and $\bm\gamma_1=\bmu$. According to \citet{Kent1982}, the contours of $f_{\rm K}$ are ``approximately ellipses'' centered at $\bm\gamma_1$, with major and minor axes $\bm\gamma_2$ and $\bm\gamma_3$, respectively. The parameter $\kappa$ controls the concentration, while $\beta$ determines the ``ovalness" of the density contours.

In the intrinsic setting, \cite{Hauberg2018} introduced an anisotropic spherical normal distribution on $\mathbb S^p$. This approach is based on the fact that any point $\bx\in\mathbb S^p\setminus\{-\bmu\}$ can be mapped to the tangent space $\mathcal T_\bmu=\{\bx\in\mathbb R^{p+1}\colon\bx^\top\bmu=0\}$ via the so-called logarithm map, which is defined as
${\rm Log}_\bmu(\bx)=\{\bx-(\bx^\top\bmu)\bmu\}d_G(\bx,\bmu)/\sin\{d_G(\bx,\bmu)\}$.
This map is an isometry, so that $d_G(\bx,\bmu)^2={\rm Log}_\bmu(\bx)^\top{\rm Log}_\bmu(\bx)$. The proposal by \citet{Hauberg2018} replaces the latter scalar product with a Mahalanobis distance in the tangent space, resulting in a distribution with density 
$f_{\rm H}(\bx;\bmu,{\bf\Lambda})=C_{\rm H}({\bf\Lambda})^{-1}\exp\{-\frac12{\Log}_\bmu(\bx)^\top{\bf\Lambda}{\rm Log}_\bmu(\bx)\},$ where $\bf\Lambda$ is a symmetric, positive semi-definite $(p+1)\times (p+1)$ matrix, having $\bmu$ as an eigenvector with associated null eigenvalue. While an exact expression for the normalizing constant $C_{\rm H}(\bf\Lambda)$ is not available, an accurate approximation was derived in \citet{Hauberg2018}.

In the following, we provide an alternative representation for \citeauthor{Hauberg2018}'s anisotropic geodesic normal distribution that does not rely on the tangent space for its formulation. We focus on  $p=2$, which will be the primary case discussed in Section \ref{sec:sph_ell}. However, most of the arguments can be easily extended to an arbitrary dimension $p$ with only minor adjustments.

\section{An alternative formulation of the anisotropic geodesic normal distribution}
\label{sec:gn}
As noted above, the matrix $\bf\Lambda$ is positive semi-definite, with a null eigenvalue corresponding to the eigenvector $\bmu$. Let us denote by $\lambda_1$ and $\lambda_2$ the remaining nonnegative eigenvalues of $\bf\Lambda$, and choose associated eigenvectors $\bm\eta$ and $\bm\xi$ so that ${\bf \Gamma}=(\bmu|\bm\eta|\bm\xi)$ is an orthogonal matrix. Thus, $\{\bmu,\bm\eta,\bm\xi\}$ is an orthonormal basis, and any $\bx\in\mathbb R^3$ can be expressed as $\bx=(\bx^\top\bmu)\bmu+(\bx^\top\bm\eta)\bm\eta+(\bx^\top\bm\xi)\bm\xi$, leading to $\bx^\top{\bf\Lambda}\bx=\lambda_1(\bx^\top\bm\eta)^2+\lambda_2(\bx^\top\bm\xi)^2$. On the other hand, for $\bx\in\mathbb S^2$ we have $d_G(\bx,\bmu)=\arccos(\bx^\top\bmu)$, resulting in $\sin^2\{d_G(\bx,\bmu)\}=1-(\bx^\top\bmu)^2=(\bx^\top\bm\eta)^2+(\bx^\top\bm\xi)^2$. Therefore, 
\begin{equation}
	\label{eq:quadratic}
	{\rm Log}_\bmu(\bx)^\top{\bm\Lambda}{\rm Log}_\bmu(\bx)=\arccos^2(\bx^\top\bmu)\frac{\lambda_1(\bx^\top\bm\eta)^2+\lambda_2(\bx^\top\bm\xi)^2}{(\bx^\top\bm\eta)^2+(\bx^\top\bm\xi)^2}.
\end{equation}
Consequently, Hauberg anisotropic geodesic normal distribution on the sphere can be alternatively defined in terms of the density
\begin{align}
	\label{eq:fgn}
	f_{\rm GN}(\bx;\lambda_1,\lambda_2,{\bf \Gamma})=C_{\rm GN}(\lambda_1,\lambda_2)^{-1}\exp\Big\{-\frac12\arccos^2(\bx^\top\bmu)\frac{\lambda_1(\bx^\top\bm\eta)^2+\lambda_2(\bx^\top\bm\xi)^2}{(\bx^\top\bm\eta)^2+(\bx^\top\bm\xi)^2}\Big\}.
\end{align}

Several important observations can be made from equation \eqref{eq:fgn}. First, it enables the definition of the geodesic normal distribution on the sphere $\mathbb S^2$ without resorting to its tangent space. Second, expression \eqref{eq:fgn} facilitates the comparison with  \citeauthor{Kent1982} distribution, since the elements needed to define the distribution are very similar, namely an orthogonal matrix $\bf\Gamma=(\bmu|\bm\eta|\bm\xi)$ and two nonnegative parameters $\lambda_1,\lambda_2$, which together account for a total of five free parameters. The vector $\bmu$ represents the mean direction, and as we will show in Section \ref{sec:sph_ell}, the distribution $f_{\rm GN}$ has true elliptical contours on the sphere (not just approximate ones), with $\bm \eta$ and $\bm\xi$ serving as the axes, and $\lambda_1,\lambda_2$ being inversely related to their lengths. In this regard, the parameters of this geodesic normal distribution have a more intuitive interpretation than those of Kent distribution. Finally, note that setting $\lambda_1=\lambda_2=\kappa$ results in the isotropic geodesic normal distribution.

\section{Spherical ellipses as density contours}
\label{sec:sph_ell}

The terms ``oval'' or ``elliptical'' are often used to describe the contours of many distributions on $\mathbb{S}^2$, although usually without a precise meaning. However, the notion of spherical ellipse has an accurate definition within the field of \emph{Spherical Geometry}, and it is a topic that has been extensively studied over the years. A recent overview is provided by \citet[][Chapter 10.1]{glaeser2016universe}.  In this section, we derive several equivalent equations for defining a spherical ellipse, which may be of mathematical interest. Furthermore, we demonstrate that the density contours of the anisotropic geodesic normal distribution indeed form true spherical ellipses. 

A spherical ellipse is defined as the locus of points on the sphere for which the sum of the geodesic distances to two fixed (non-antipodal) points remains constant. Specifically, given two points $\bbf_1,\bbf_2\in\mathbb S^2$ with geodesic distance $d_G(\bbf_1,\bbf_2)=\arccos({\bbf_1}\!\!^\top\bbf_2)=2\gamma<\pi$, and an angle $\alpha$, the spherical ellipse with foci $\bbf_1,\bbf_2$ and major axis length $2\alpha$ is defined as the set of points $\bx\in\mathbb S^2$ satisfying
\begin{equation}\label{eq:se1}
	d_G(\bx,\bbf_1)+d_G(\bx,\bbf_2)=2\alpha.
\end{equation}

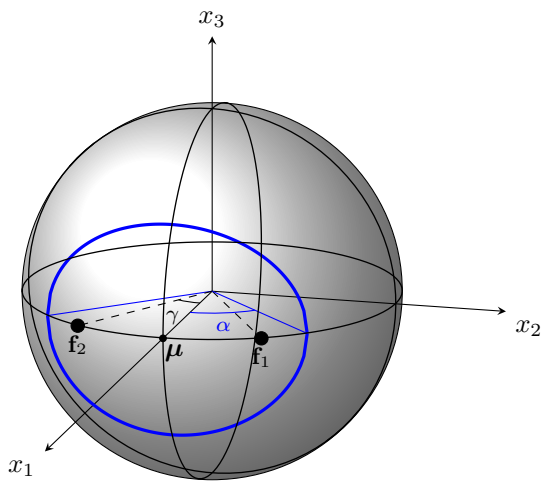
\begin{figure}
	\centering
	\begin{tikzpicture}[scale=2.5,tdplot_main_coords]
		
		\shadedraw[tdplot_screen_coords,ball color = white] (0,0) circle (1);
		
		\draw[-stealth] (0,0,0) -- (3.4,0,0)
		node[below left] {$x_1$};
		\draw[-stealth] (0,0,0) -- (0,1.60,0)
		node[below right] {$x_2$};
		\draw[-stealth] (0,0,0) -- (0,0,1.40)
		node[above] {$x_3$};
		
		\coordinate (O) at (0,0,0);
		
		\coordinate (F) at ({sqrt(3)/2}, 0.5, 0);
		\draw[dashed] (O) --  (F)
		node[below] {\small $\bbf_1$};
		\draw[fill=black] (F) circle (1pt);
		
		\tdplotdrawarc{(O)}{0.25}{0}{-30}{anchor=north east}{\footnotesize$\gamma$};
		
		\coordinate (A) at ({sqrt(2)/2}, {sqrt(2)/2}, 0);
		\draw[thin,color=blue] (O) --  (A);
		\tdplotdrawarc[color=blue]{(O)}{0.45}{0}{45}{anchor=north}{\footnotesize$\alpha$};
		\draw[thin,color=blue] (O) --  ({sqrt(2)/2}, {-sqrt(2)/2}, 0);
		
		\coordinate (G) at ({sqrt(3)/2}, -0.5, 0);
		\draw[dashed] (O) -- (G)
		node[below] {\small $\bbf_2$};
		\draw[fill=black] (G) circle (1pt);
		
		\draw[fill=black] (1,0,0) circle (0.5pt) node[below, xshift=.15cm, yshift=.05cm] {$\bmu$};
		
		
		\draw[blue, domain=-45:45, variable=\t, samples=101, very thick] plot
		({cos(\t)*sqrt(2/(1+2*cos(\t)*cos(\t)))},
		{sin(\t)*sqrt(2/(1+2*cos(\t)*cos(\t)))},
		{sqrt(1-2/(1+2*cos(\t)*cos(\t)))});
		
		\draw[blue, domain=-45:45, variable=\t, samples=101, very thick] plot
		({cos(\t)*sqrt(2/(1+2*cos(\t)*cos(\t)))},
		{sin(\t)*sqrt(2/(1+2*cos(\t)*cos(\t)))},
		{-sqrt(1-2/(1+2*cos(\t)*cos(\t)))});

		\pgfmathsetmacro{\CameraX}{sin(\ViewAzimuth)*cos(\ViewElevation)}
		\pgfmathsetmacro{\CameraY}{-cos(\ViewAzimuth)*cos(\ViewElevation)}
		\pgfmathsetmacro{\CameraZ}{sin(\ViewElevation)}
		\path[use as bounding box] (-1,-1) rectangle (1,1); 
		
		\begin{axis}[
			hide axis,
			view={\ViewAzimuth}{\ViewElevation},     
			every axis plot/.style={very thin},
			disabledatascaling,                      
			anchor=origin,                           
			viewport={\ViewAzimuth}{\ViewElevation}, 
			]
			\addFGBGplot[domain=0:2*pi, samples=100, samples y=1] ({cos(deg(x))}, {sin(deg(x))}, 0);
			\addFGBGplot[domain=0:2*pi, samples=100, samples y=1] (0, {sin(deg(x))}, {cos(deg(x))});
			\addFGBGplot[domain=0:2*pi, samples=100, samples y=1] ({sin(deg(x))}, 0, {cos(deg(x))});
		\end{axis}
	\end{tikzpicture}
	\caption{Spherical ellipse with foci $\bbf_1=(\cos\gamma,\sin\gamma,0)$ and $\bbf_2=(\cos\gamma,-\sin\gamma,0)$ for $\gamma=\pi/6$, center $\bmu=(1,0,0)$ and $\alpha=\pi/4$.}
	\label{fig:sph_ell}
\end{figure}

The great circle on the sphere passing through the two foci $\bbf_1$ and $\bbf_2$ determines the major axis of the ellipse. The center of the ellipse is the midpoint of the shortest arc joining $\bbf_1$ and $\bbf_2$, that is, $\bmu=(\bbf_1+\bbf_2)/\|\bbf_1+\bbf_2\|=(\bbf_1+\bbf_2)/(2\cos\gamma)\in\mathbb S^2$. Thus, the angle $\gamma$ represents the geodesic distance from either focus to the ellipse center, and must satisfy $\gamma<\alpha<\pi-\gamma$. Figure~\ref{fig:sph_ell} shows a spherical ellipse with foci given by $\bbf_1=(\cos\gamma,\sin\gamma,0)$ and $\bbf_2=(\cos\gamma,-\sin\gamma,0)$ for $\gamma=\pi/6$, center $\bmu=(1,0,0)$ and $\alpha=\pi/4$.

Although \eqref{eq:se1} provides the most straightforward definition of a spherical ellipse, alternative forms are often more convenient. Instead of directly using the foci, the ellipse can be characterized by its center $\bmu$ and the direction $\bm\eta$ defining the major axis. To determine $\bm\eta$, we need to find a unit vector from the space spanned by $\bbf_1$ and $\bbf_2$ that is orthogonal to $\bmu$. It can be checked that the vector $\bm\eta=(\bbf_1-\bbf_2)/(2\sin\gamma)$ satisfies the required conditions. Next, express the foci $\bbf_1$ and $\bbf_2$ in terms of $\bmu$ and $\bm\eta$, resulting in $\bbf_1=\bm\mu\cos\gamma+\bm\eta\sin \gamma$ and $\bbf_2=\bm\mu\cos\gamma-\bm\eta\sin \gamma$. Substituting these into the original ellipse \eqref{eq:se1}, we obtain
\begin{equation}\label{eq:se2}
	\arccos\big(\bx^\top\bm\mu\cos\gamma+\bx^\top\bm\eta\sin\gamma\big)+\arccos\big(\bx^\top\bm\mu\cos\gamma-\bx^\top\bm\eta\sin\gamma\big)=2\alpha,
\end{equation}
with the constraints $0\leq\gamma<\pi/2$ and $\gamma<\alpha<\pi-\gamma$, where $\bm\eta\in\mathbb S^2$ is orthogonal to $\bmu$. In the case where $\gamma=0$, this equation simplifies to a spherical circle with center $\bmu$ and radius $\alpha<\pi$, as it reduces to $\arccos(\bx^\top\bm\mu)=\alpha$.

Alternatively, the ellipse equation can be written in terms of the semi-minor axis length $\beta$. Using the spherical Pythagorean theorem for the minor axis vertex it follows that $\cos \alpha = \cos \beta \cos \gamma$, so we find that
$\cos\gamma=(\cos\alpha)/(\cos\beta)$ and $\sin\gamma=\{1-(\cos^2\alpha)/(\cos^2\beta)\}^{1/2}$. Substituting these expressions into equation (\ref{eq:se2}) gives
\begin{align}\label{eq:se3}
	\arccos\big\{\bx^\top\bm\mu\tfrac{\cos\alpha}{\cos\beta}+\bx^\top\bm\eta(1-\tfrac{\cos^2\alpha}{\cos^2\beta})^{1/2}\big\}+\arccos\big\{\bx^\top\bm\mu\tfrac{\cos\alpha}{\cos\beta}-\bx^\top\bm\eta(1-\tfrac{\cos^2\alpha}{\cos^2\beta})^{1/2}\big\}=2\alpha,
\end{align}
with the constraints $0<\beta\leq\alpha<\pi$, and $\bm\eta\in\mathbb S^2$ orthogonal to $\bmu\in\mathbb S^2$. Note that this equation involves five free parameters, similar to the general equation of an ellipse in $\mathbb R^2$. In this case, the equation of a spherical circle with radius $\rho$ is obtained by setting $\alpha=\beta=\rho$.

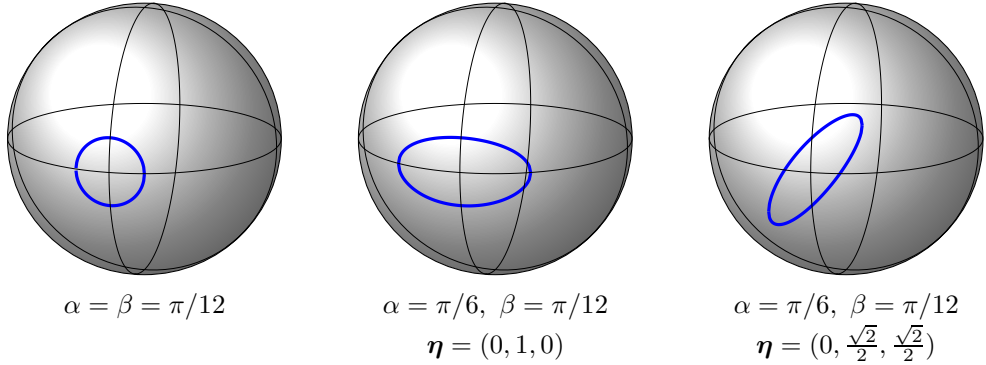
\begin{figure}
	\centering
	\begin{tabular}{@{}c@{\hskip 1cm}c@{\hskip 1cm}c@{}}

		\begin{tikzpicture}[scale=1.8,tdplot_main_coords]
			
			\shadedraw[tdplot_screen_coords,ball color = white] (0,0) circle (1);

			\coordinate (O) at (0,0,0);
			
			\coordinate (F) at ({sqrt(3)/2}, 0.5, 0);
			
			\coordinate (A) at ({sqrt(2)/2}, {sqrt(2)/2}, 0);

			\coordinate (G) at ({sqrt(3)/2}, -0.5, 0);

			
			\draw[blue, domain=-15:15, variable=\t, samples=101, very thick] plot
			({cos(15)},
			{cos(15)*tan(\t)},
			{sqrt(1-cos(15)^2/cos(\t)^2});
			
			\draw[blue, domain=-15:15, variable=\t, samples=101, very thick] plot
			({cos(15)},
			{cos(15)*tan(\t)},
			{-sqrt(1-cos(15)^2/cos(\t)^2});

			\pgfmathsetmacro{\CameraX}{sin(\ViewAzimuth)*cos(\ViewElevation)}
			\pgfmathsetmacro{\CameraY}{-cos(\ViewAzimuth)*cos(\ViewElevation)}
			\pgfmathsetmacro{\CameraZ}{sin(\ViewElevation)}
			\path[use as bounding box] (-1,-1) rectangle (1,1); 
			
			\begin{axis}[
				hide axis,
				view={\ViewAzimuth}{\ViewElevation},     
				every axis plot/.style={very thin},
				disabledatascaling,                      
				anchor=origin,                           
				viewport={\ViewAzimuth}{\ViewElevation}, 
				]
				\addFGBGplot[domain=0:2*pi, samples=100, samples y=1] ({cos(deg(x))}, {sin(deg(x))}, 0);
				\addFGBGplot[domain=0:2*pi, samples=100, samples y=1] (0, {sin(deg(x))}, {cos(deg(x))});
				\addFGBGplot[domain=0:2*pi, samples=100, samples y=1] ({sin(deg(x))}, 0, {cos(deg(x))});
			\end{axis}
		\end{tikzpicture}&
		
		\begin{tikzpicture}[scale=1.8,tdplot_main_coords]
			
			\shadedraw[tdplot_screen_coords,ball color = white] (0,0) circle (1);

			
			\draw[blue, domain=-30:30, variable=\t, samples=101, very thick] plot
			({cos(\t)*sqrt(-(2+sqrt(3))/(-2*sqrt(3)+2*(-2+sqrt(3))*cos(2*\t)))},
			{sin(\t)*sqrt(-(2+sqrt(3))/(-2*sqrt(3)+2*(-2+sqrt(3))*cos(2*\t)))},
			{sqrt(1+(2+sqrt(3))/(-2*sqrt(3)+2*(-2+sqrt(3))*cos(2*\t)))});
			
			\draw[blue, domain=-30:30, variable=\t, samples=101, very thick] plot
			({cos(\t)*sqrt(-(2+sqrt(3))/(-2*sqrt(3)+2*(-2+sqrt(3))*cos(2*\t)))},
			{sin(\t)*sqrt(-(2+sqrt(3))/(-2*sqrt(3)+2*(-2+sqrt(3))*cos(2*\t)))},
			{-sqrt(1+(2+sqrt(3))/(-2*sqrt(3)+2*(-2+sqrt(3))*cos(2*\t)))});

			\pgfmathsetmacro{\CameraX}{sin(\ViewAzimuth)*cos(\ViewElevation)}
			\pgfmathsetmacro{\CameraY}{-cos(\ViewAzimuth)*cos(\ViewElevation)}
			\pgfmathsetmacro{\CameraZ}{sin(\ViewElevation)}
			\path[use as bounding box] (-1,-1) rectangle (1,1); 
			
			\begin{axis}[
				hide axis,
				view={\ViewAzimuth}{\ViewElevation},     
				every axis plot/.style={very thin},
				disabledatascaling,                      
				anchor=origin,                           
				viewport={\ViewAzimuth}{\ViewElevation}, 
				]
				\addFGBGplot[domain=0:2*pi, samples=100, samples y=1] ({cos(deg(x))}, {sin(deg(x))}, 0);
				\addFGBGplot[domain=0:2*pi, samples=100, samples y=1] (0, {sin(deg(x))}, {cos(deg(x))});
				\addFGBGplot[domain=0:2*pi, samples=100, samples y=1] ({sin(deg(x))}, 0, {cos(deg(x))});
			\end{axis}
		\end{tikzpicture}&

		\begin{tikzpicture}[scale=1.8,tdplot_main_coords]
			
			\shadedraw[tdplot_screen_coords,ball color = white] (0,0) circle (1);

			
			\draw[blue, domain=-30:30, variable=\t, samples=101, very thick] plot
			({cos(\t)*sqrt(-(2+sqrt(3))/(-2*sqrt(3)+2*(-2+sqrt(3))*cos(2*\t)))},
			{sin(\t)*sqrt(-(2+sqrt(3))/(-2*sqrt(3)+2*(-2+sqrt(3))*cos(2*\t)))*sqrt(2)/2},
			{sin(\t)*sqrt(-(2+sqrt(3))/(-2*sqrt(3)+2*(-2+sqrt(3))*cos(2*\t)))*sqrt(2)/2+sqrt(1+(2+sqrt(3))/(-2*sqrt(3)+2*(-2+sqrt(3))*cos(2*\t)))});
			
			\draw[blue, domain=-30:30, variable=\t, samples=101, very thick] plot
			({cos(\t)*sqrt(-(2+sqrt(3))/(-2*sqrt(3)+2*(-2+sqrt(3))*cos(2*\t)))},
			{sin(\t)*sqrt(-(2+sqrt(3))/(-2*sqrt(3)+2*(-2+sqrt(3))*cos(2*\t)))*sqrt(2)/2},
			{sin(\t)*sqrt(-(2+sqrt(3))/(-2*sqrt(3)+2*(-2+sqrt(3))*cos(2*\t)))*sqrt(2)/2-sqrt(1+(2+sqrt(3))/(-2*sqrt(3)+2*(-2+sqrt(3))*cos(2*\t)))});

			\pgfmathsetmacro{\CameraX}{sin(\ViewAzimuth)*cos(\ViewElevation)}
			\pgfmathsetmacro{\CameraY}{-cos(\ViewAzimuth)*cos(\ViewElevation)}
			\pgfmathsetmacro{\CameraZ}{sin(\ViewElevation)}
			\path[use as bounding box] (-1,-1) rectangle (1,1); 
			
			\begin{axis}[
				hide axis,
				view={\ViewAzimuth}{\ViewElevation},     
				every axis plot/.style={very thin},
				disabledatascaling,                      
				anchor=origin,                           
				viewport={\ViewAzimuth}{\ViewElevation}, 
				]
				\addFGBGplot[domain=0:2*pi, samples=100, samples y=1] ({cos(deg(x))}, {sin(deg(x))}, 0);
				\addFGBGplot[domain=0:2*pi, samples=100, samples y=1] (0, {sin(deg(x))}, {cos(deg(x))});
				\addFGBGplot[domain=0:2*pi, samples=100, samples y=1] ({sin(deg(x))}, 0, {cos(deg(x))});
			\end{axis}
		\end{tikzpicture}\\

		$\alpha=\beta=\pi/12$ & $\alpha=\pi/6,\ \beta=\pi/12$ & $\alpha=\pi/6,\ \beta=\pi/12$ \\
		& $\bm\eta=(0,1,0)$ & $\bm\eta=(0,\frac{\sqrt2}2,\frac{\sqrt2}2)$

	\end{tabular}
	\caption{
		Spherical ellipses centered at $\bmu=(1,0,0)$ with different parameter configurations.
	}
	\label{fig:sph_ell2}
\end{figure}

Figure~\ref{fig:sph_ell2} shows three examples of spherical ellipses centered at $\bmu=(1,0,0)$, each with different parameter configurations. The first one (left panel), with $\alpha=\beta=\pi/12$, corresponds to a spherical circle. The second one (center panel), with $\bm\eta=(0,1,0)$ and $\alpha=\pi/6$, $\beta=\pi/12$, shows an ellipse aligned with the equator and the reference meridian. The third one (right panel) has the same axis lengths but a different orientation for the major axis, with $\bm\eta=(0,\frac{\sqrt2}2,\frac{\sqrt2}2)$.

Some references mention even simpler equations in case the ellipse is assumed to be in ``standard position''. In $\mathbb R^2$, this refers to an ellipse centered at the origin and aligned with the coordinate axes. In the spherical case, the standard ellipse can be considered to be centered at $\bmu=(1,0,0)$, with the major axis along the direction of $\bm\eta=(0,1,0)$. Then, writing $\bx=(x_1,x_2,x_3)\in\mathbb R^3$, equation \eqref{eq:se2} simplifies to $\arccos(x_1\cos\gamma+x_2\sin\gamma)+\arccos(x_1\cos\gamma-x_2\sin\gamma)=2\alpha$.
Next, after some straightforward calculations we get
\begin{equation}\label{eq:se4}
	x_1^2\frac{\cos^2\gamma}{\cos^2\alpha}+x_2^2\frac{\sin^2\gamma}{\sin^2\alpha}=1,
\end{equation}
as shown in \citet[][pp. 446--448]{glaeser2016universe}. This recalls the usual equation $x_1^2/a^2+x_2^2/b^2=1$ of a standard ellipse in $\mathbb R^2$.

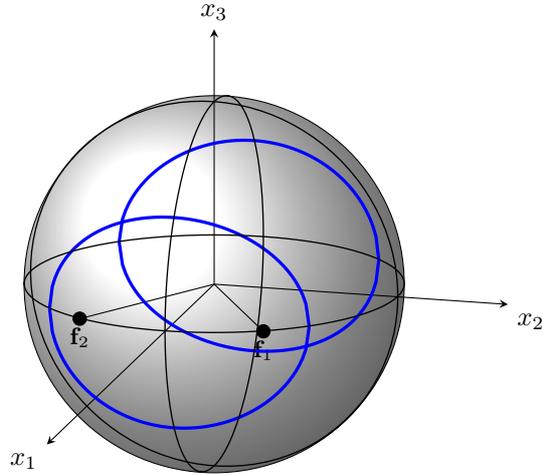
\begin{figure}
	\centering
	\begin{tikzpicture}[scale=2.5,tdplot_main_coords]
		
		\shadedraw[tdplot_screen_coords,ball color = white] (0,0) circle (1);
		
		\draw[-stealth] (0,0,0) -- (3.4,0,0)
		node[below left] {$x_1$};
		\draw[-stealth] (0,0,0) -- (0,1.60,0)
		node[below right] {$x_2$};
		\draw[-stealth] (0,0,0) -- (0,0,1.40)
		node[above] {$x_3$};
		
		\coordinate (O) at (0,0,0);
		
		\coordinate (F) at ({sqrt(3)/2}, 0.5, 0);
		\draw[thin] (O) --  (F)
		node[below] {\small $\bbf_1$};
		\draw[fill=black] (F) circle (1pt);
		
		
		\coordinate (G) at ({sqrt(3)/2}, -0.5, 0);
		\draw[thin] (O) -- (G)
		node[below] {\small $\bbf_2$};
		\draw[fill=black] (G) circle (1pt);
		
		
		\draw[blue, domain=-45:45, variable=\t, samples=101, very thick] plot
		({cos(\t)*sqrt(2/(1+2*cos(\t)*cos(\t)))},
		{sin(\t)*sqrt(2/(1+2*cos(\t)*cos(\t)))},
		{sqrt(1-2/(1+2*cos(\t)*cos(\t)))});
		
		\draw[blue, domain=-45:45, variable=\t, samples=101, very thick] plot
		({cos(\t)*sqrt(2/(1+2*cos(\t)*cos(\t)))},
		{sin(\t)*sqrt(2/(1+2*cos(\t)*cos(\t)))},
		{-sqrt(1-2/(1+2*cos(\t)*cos(\t)))});

		
		\draw[blue, domain=-45:45, variable=\t, samples=101, very thick, opacity=0.25] plot
		(-{cos(\t)*sqrt(2/(1+2*cos(\t)*cos(\t)))},
		{sin(\t)*sqrt(2/(1+2*cos(\t)*cos(\t)))},
		{sqrt(1-2/(1+2*cos(\t)*cos(\t)))});
		
		\draw[blue, domain=-45:45, variable=\t, samples=101, very thick, opacity=0.25] plot
		(-{cos(\t)*sqrt(2/(1+2*cos(\t)*cos(\t)))},
		{sin(\t)*sqrt(2/(1+2*cos(\t)*cos(\t)))},
		{-sqrt(1-2/(1+2*cos(\t)*cos(\t)))});

		\pgfmathsetmacro{\CameraX}{sin(\ViewAzimuth)*cos(\ViewElevation)}
		\pgfmathsetmacro{\CameraY}{-cos(\ViewAzimuth)*cos(\ViewElevation)}
		\pgfmathsetmacro{\CameraZ}{sin(\ViewElevation)}
		\path[use as bounding box] (-1,-1) rectangle (1,1); 
		
		\begin{axis}[
			hide axis,
			view={\ViewAzimuth}{\ViewElevation},     
			every axis plot/.style={very thin},
			disabledatascaling,                      
			anchor=origin,                           
			viewport={\ViewAzimuth}{\ViewElevation}, 
			]
			\addFGBGplot[domain=0:2*pi, samples=100, samples y=1] ({cos(deg(x))}, {sin(deg(x))}, 0);
			\addFGBGplot[domain=0:2*pi, samples=100, samples y=1] (0, {sin(deg(x))}, {cos(deg(x))});
			\addFGBGplot[domain=0:2*pi, samples=100, samples y=1] ({sin(deg(x))}, 0, {cos(deg(x))});
		\end{axis}
	\end{tikzpicture}
	\caption{Graphical representation of equation \eqref{eq:se4} using the same parameter configuration as in Figure~\ref{fig:sph_ell}, resulting in two antipodal ellipses.}
	\label{fig:two_ell}
\end{figure}

However, in contrast to the Euclidean case, the transition from equation \eqref{eq:se3} to \eqref{eq:se4} introduces a spurious curve. While the set of points satisfying equation \eqref{eq:se3} corresponds to a single curve on the sphere (as shown previously in Figure~\ref{fig:sph_ell}), the locus of points satisfying equation \eqref{eq:se4} consists of two antipodal ellipses, centered at $\bmu$ and $-\bmu$, respectively, as shown in Figure~\ref{fig:two_ell}. This outcome is not so surprising, since equation \eqref{eq:se4} corresponds to a elliptic cylinder in $\mathbb R^3$, and its intersection with $\mathbb S^2$ gives rise to the two spherical ellipses.

To address the previous issue, we derive an alternative equation that includes only one spherical ellipse. Our development is based on the tangent space $\mathcal T_\bmu$, but no mapping to $\mathcal T_\bmu$ is eventually required for the final equation. As discussed in Section \ref{sec:normal}, let $\bm\eta,\bm\xi\in\mathbb S^2$ be vectors such that $\{\bmu,\bm\eta,\bm\xi\}$ form an orthonormal basis of $\mathbb R^3$ and consider the orthogonal matrix ${\bf\Gamma}=(\bmu|\bm\eta|\bm\xi)$. Let $\mathcal E\subset \mathcal S^2$ represent the locus of points $\bx\in\mathbb S^2$ satisfying
\begin{equation}\label{eq:se5}
	\arccos^2(\bx^\top\bmu)\frac{(\bx^\top\bm\eta)^2/\alpha^2+(\bx^\top\bm\xi)^2/\beta^2}{(\bx^\top\bm\eta)^2+(\bx^\top\bm\xi)^2}=1.
\end{equation}
In view of \eqref{eq:quadratic}, the image of $\mathcal E$ through $\Log_\bmu$ forms an ellipse $\mathcal E^*$ in the tangent space $\mathcal T_\bmu$, since it consists of all the points $\by=\Log_\bmu(\bx)\in\mathcal T_\bmu$ satisfying $\by^\top{\bf\Lambda}\by=1$, where ${\bf\Lambda}={\bf \Gamma}{\bf D}{\bf\Gamma}^\top$ and $\mathbf D$ is the diagonal matrix ${\rm diag}(0,1/\alpha^2,1/\beta^2)$. That is, in terms of the basis $\{\bm\eta,\bm\xi\}$ of the vector space $\mathcal T_\bmu$, the set $\mathcal E^*$ is the ellipse with equation $(\by^\top\bm\eta)^2/\alpha^2+(\by^\top\bm\xi)^2/\beta^2=1$, where the major and minor axes are determined by $\bm\eta$ and $\bm\xi$, with lengths $\alpha$ and $\beta$, respectively.
Then, the foci of the ellipse $\mathcal E^*\subset\mathcal T_\bmu$ are located at $\gamma\bm\eta$ and $-\gamma\bm\eta$, with $0\leq\gamma<\alpha$. Therefore,  $\mathcal E^*$ can be equivalently described as the locus of points $\by\in\mathcal T_\bmu$ such that $d_\mathcal T(\by,\gamma\bm\eta)+d_\mathcal T(\by,-\gamma\bm\eta)=2\alpha$, where $d_\mathcal T$ denotes the Euclidean distance in $\mathcal T_\bmu$, defined as $d_\mathcal T(\by_1,\by_2)^2=(\by_1-\by_2)^\top(\by_1-\by_2)$ for $\by_1,\by_2\in\mathcal T_\bmu$. Since $\Log_\bmu$ is an isometry, this implies that $\mathcal E$ is the locus of points $\bx\in\mathbb S^2$ such that $d_G(\bx,\bbf_1)+d_G(\by,\bbf_2)=2\alpha$, where $\bbf_1=\Log_\bmu^{-1}(\gamma\bm\eta)$ and $\bbf_2=\Log_\bmu^{-1}(-\gamma\bm\eta)$. Hence, $\mathcal E$ is a spherical ellipse.

Thus, \eqref{eq:se5} represents the equation of the spherical ellipse, with axes 
determined by $\bm\eta$ and $\bm\xi$ and semi-axis lengths $\alpha$ and $\beta$, respectively. This demonstrates, on the one hand, that the density contours of the geodesic normal distribution are spherical ellipses. On the other hand, equation \eqref{eq:se5} offers an alternative equation for a spherical ellipse, which appears to be previously unknown.

\bigskip

\noindent{\bf Acknowledgments}

Both authors are supported by grant PID2021-124051NB-I00 from the Spanish Ministerio de Ciencia e Innovación. The first author also acknowledges support by the grant PID2023-148081NB-I00 from the Spanish Ministerio de Ciencia, Innovación y Universidades.

\bibliographystyle{apalike}

\end{document}